\documentclass{article}
\usepackage{amssymb,amsmath}
\usepackage{graphics}
\input xy
\xyoption{all}
\usepackage[all,cmtip]{xy}

\begin{document}

\title{Elliptic complexes over $C^*$-algebras of compact 
operators}
\author{Svatopluk  Kr\'ysl \footnote{{\it E-mail address}: 
Svatopluk.Krysl@mff.cuni.cz}\\ {\it \small Faculty of Mathematics and Physics, 
Charles University in Prague, Czech Republic}
\thanks{The author thanks for a financial support from the foundation PRVOUK at 
the School of Mathematics of the Faculty of Mathematics and Physics 
granted by the Charles University in Prague.}}

\maketitle \noindent

\centerline{\large\bf Abstract}  

For a $C^*$-algebra $A$ of compact operators  and  a compact manifold $M,$
we prove that  the Hodge theory 
holds for
$A$-elliptic complexes of pseudodifferential operators acting on smooth sections
of finitely generated projective $A$-Hilbert bundles over $M.$ 
For these $C^*$-algebras, we get also a 
topological isomorphism between the cohomology groups 
of an $A$-elliptic complex 
and the
space of harmonic elements. Consequently,   the cohomology groups appear 
to be finitely generated projective $C^*$-Hilbert modules and especially,
Banach spaces. 
We prove as well, that if the Hodge theory holds for a complex in the category 
of Hilbert $A$-modules and continuous adjointable Hilbert $A$-module 
homomorphisms, 
the complex is self-adjoint parametrix possessing.

\bigskip 
 
{\it Math.  Subj.  Class. 2010:} Primary 46M18; Secondary 46L08, 46M20

{\it Key words:} Hodge theory, Hilbert $C^*$-modules, $C^*$-Hilbert bundles, 
elliptic
systems of partial differential equation

\section{Introduction}

This paper is a continuation of  papers 
\cite{KryslAGAG1} and \cite{KryslAGAG2}, devoted to the Hodge theory for 
Hilbert and pre-Hilbert $C^*$-modules and to an  application of this theory
to $A$-elliptic complexes of operators acting on sections of specific  
$C^*$-Hilbert 
bundles over compact manifolds.

Let $A$ be a $C^*$-algebra and $M$ be a compact manifold. In 
\cite{KryslAGAG2}, the Hodge theory is proved to hold for an arbitrary 
 $A$-elliptic complex of operators acting on smooth sections of  finitely 
generated projective 
$A$-Hilbert bundles over $M$ if
 the images of the extensions to the  Sobolev spaces  of the 
Laplacians of 
the complex
 are closed.
One of the  main results achieved in this paper is that one can omit
 the assumption on the images if $A$ is a
$C^*$-algebra of compact operators and still get the claim of the theorem.
 
We define what it means that {\it the
Hodge theory holds} for a complex in an additive and dagger 
category and
study this concept in a detail in categories of 
pre-Hilbert and Hilbert modules and continuous adjointable $A$-equivariant   
maps. These categories constitute a  special class of the 
so-called 
$R$-module categories that are in addition, equipped with an involution on the 
morphisms spaces.
 We  say that the {\it Hodge theory holds} for a complex $d^{\bullet}=(U^i, 
d_i: 
U^i \to 
U^{i+1})_{i\in \mathbb{Z}}$
  in an additive and dagger category $\mathfrak{C}$ or that $d^{\bullet}$ is of 
{\it Hodge type} if  for each
	$i \in \mathbb{Z},$ we have $$U^i = \textrm{Im}\ d_{i-1} \oplus 
\textrm{Im}\, d_i^* \oplus 
\, \textrm{Ker} 
\, \Delta_i,$$ 
where $\Delta_i = d_i^* d_i + d_{i-1} d_{i-1}^*,$ and $d_i^*$ and 
$d_{i-1}^*$ 
are the adjoints
of  $d_i$ and $d_{i-1},$ respectively. 
The operators $\Delta_i,$ $i\in 
\mathbb{Z},$ are called the 
{\it Laplace operators} of $d^{\bullet}.$ 
The term  "dagger category" is explained in the paper.

  For a  $C^*$-algebra $A$,  we  consider the   
category 
$PH_A^*$ of right pre-Hilbert $A$-modules and continuous adjointable
$A$-equivariant maps.
 The full subcategory of $PH_A^*,$ the object of which are right Hilbert 
$A$-modules is denoted  by $H_A^*$ and it is called the category of Hilbert 
$A$-modules. 
See Kaplansky \cite{Kaplansky}, Paschke \cite{Paschke}, Lance \cite{Lance} and 
Manuilov, Troitsky \cite{MT} 
for information on (pre-)Hilbert modules.
Recall that each object in $PH_A^*$ inherits a norm derived from the 
$A$-product 
defined on it. The categories $PH_A^*$ and $H_A^*$ are 
additive and dagger with respect to the orthogonal direct sum and an
involution defined by the $A$-product.
 
In Kr\'ysl \cite{KryslAGAG2}, the  so-called self-adjoint parametrix possessing
complexes in $PH_A^*$ are introduced.
According to results in  that paper, any self-adjoint parametrix 
possessing complex in $PH_A^*$  is of Hodge type  and its cohomology groups are 
pre-Hilbert $A$-modules  isomorphic to the kernels of the Laplace 
ope-\\rators as pre-Hilbert  
$A$-modules.  Especially, the cohomology groups are normed spaces. 
In the present paper, we 
prove the opposite  implication in the category $H_A^*$, i.e., that if 
the Hodge  theory holds for a  
complex in $H_A^*,$ the complex  is self-adjoint 
parametrix  possessing.   Thus, in $H_A^*$ the condition of  being  
self-adjoint parametrix possessing characterizes the  Hodge type complexes.

  Let us recall that the Hodge theory is well known to hold for 
elliptic complexes  of pseudodifferential operators acting 
on smooth sections of finite rank vector bundles over compact 
manifolds. Classical examples are  deRham and  Dolbeault complexes 
over compact manifolds. See, e.g., Palais \cite{Palais} or Wells \cite{Wells}.
Fomenko, Mishchenko prove in \cite{FM}, that the continuous extensions  
of an $A$-elliptic operator to the  Sobolev section spaces are 
$A$-Fredholm. In \cite{Guljas}, Baki\'c and Gulja\v{s} prove that any 
$A$-Fredholm endomorphism $F: U \to U$ in $H_A^*$ has closed image if  $A$ is a 
$C^*$-algebra of compact 
operators. 
By a simple transforming, we generalize this result to the case  
of an
$A$-Fredholm morphism  $F:U \to 
V$ acting  between  Hilbert $A$-modules $U$ and $V.$
 In particular, we prove that the 
image of $F$ is closed. For $C^*$-algebras of compact operators, we further derive a transfer theorem 
which roughly speaking, enables us to deduce certain properties of   pre-Hilbert 
$A$-module maps from the appropriate properties of their extensions. We specify the properties and the maps in the theorem formulation. 
Applying the mentioned  theorem generalizing the result of Baki\'c and 
Gulja\v{s}, 
we get that the images of the extended Laplace operators of an $A$-elliptic 
complex are closed. The transfer theorem enables us to prove that in this case, the Laplace 
operators  themselves have closed images, they are self-adjoint parametrix 
possessing, and 
consequently, that the complex is of Hodge type.

The motivation for our research comes from 
quantum field theories which  aim to include constraints -- especially, from 
the
Becchi, Rouet, Stora and Tyutin or simply BRST quantization. 
See Henneaux, Teitelboim \cite{HT}, Horuzhy, Voronin \cite{Horuzhy}, 
Carchedi, Roytenberg \cite{Roytenberg} and the references there.
Let us explain the connection in a more detail.
In the BRST quantization, one constructs complexes whose cohomology 
groups  represent state spaces of a given physical system.   Because the state 
spaces in 
quantum theories  are usually formed by infinite dimensional vector spaces, the 
 co-cycle spaces  for the cohomology groups have to be infinite dimensional as 
well. 
It is agreed that the state spaces shall be equipped  
with a topology   because of the testing of the theory by 
measurements. Since the measurements do not give the precise value of 
a measured  observable  (a result of a measurement is always a value 
together with an error estimate),
the state spaces  shall have good  a 
good behavior of limits of converging sequences.  Especially, it is 
desirable that the limit of 
a converging sequence is  unique. It is well known that the uniqueness of 
limits 
in a  topological space forces the space to be T1. However, the T1 separation 
axiom in a topological  vector space implies that the topological vector
space is already Hausdorff. (For it, see, e.g., Theorem 1.12 in Rudin 
\cite{Rudin}.)
The quotient of a topological vector space  is non-Hausdorff in the quotient 
topology
if and only  if  the space by which one divides is not closed. 
If we insist that the 
state spaces are cohomology groups, we shall be able 
to assure that the spaces of co-boundaries are closed. For a explanation  of 
the requirements on a physical theory considered
above,  we refer to Ludwig \cite{Ludwig} and to a still appealing 
paper of von Neumann \cite{Neumann}.
We hope that our work can be  relevant for physics at least in the case when 
a particular BRST complex  appears to be self-adjoint parametrix possessing  
in 
the categories  $PH_A^*$ or $H_A^*$ for an arbitrary $C^*$-algebra $A,$ or
an $A$-elliptic complex in finitely generated projective $A$-Hilbert bundles 
over a compact manifold if  $A$ is a $C^*$-algebra of compact operators.
A further inspiring topic from physics is the   
parallel transport in Hilbert bundles considered in a connection with quantum theory. See, 
e.g., Drechsler, Tuckey \cite{Drechsler}.  

Let us notice that in Troitsky \cite{Troitsky}, indices of $A$-elliptic 
complexes are investigated. In that paper 
the operators are, quite naturally, allowed to be changed by an $A$-compact 
perturbation in order the index is an element of the appropriate $K$-group. See also Schick \cite{Schick}.
If the reader is interested in a possible application of the Hodge theory for 
$A$-elliptic complexes, we 
refer to  Kr\'ysl \cite{KryslDGA}.

In the second chapter, we give a definition of the Hodge type complex, 
recall definitions of a pre-Hilbert and a Hilbert 
$C^*$-module, and give several examples of them. 
We prove that   complexes in the category of Hilbert spaces and continuous 
maps are of Hodge type if the images of their Laplace operators are closed (Lemma 1).
Further, we recall the definition of a self-adjoint parametrix 
complex in $PH_A^*$ and some of its properties including the fact that they are 
of Hodge 
type (Theorem 2). We  prove that if  a 
complex in $H_A^*$ is of Hodge type, it is already self-adjoint parametrix 
possessing 
(Theorem 3).
At the end of the second section, we give   examples of complexes the cohomology groups of
which are not Hausdorff spaces.
In the third chapter, we summarize the result of
Baki\'{c} and Gulja\v{s} (Theorem 4), give the mentioned generalization of it 
(Corollary 5), and prove the transfer theorem (Theorem 6). In the fourth 
section, basic facts on differential operators acting on sections of 
$A$-Hilbert bundles over compact manifolds are recalled. 
In this chapter, the theorem on properties of $A$-elliptic complexes in 
finitely generated projective $A$-Hilbert bundles over compact manifolds 
is  proved (Theorem 9).

\bigskip

\noindent {\bf Preamble}: All manifolds and bundles are assumed to be smooth. 
The base 
manifolds of bundles are assumed to be finite dimensional.
When an index of a labeled object  exceeds its allowed range, 
the object is set 
to be zero. We do not suppose the Hilbert spaces to be separable.  

\section{Self-adjoint maps and complexes possessing a parametrix}

Let us recall that a category $\mathfrak{C}$ is called a {\it dagger category} 
if
there is a contra\--variant functor $*: \mathfrak{C}\to \mathfrak{C}$  
which is the identity on the objects 
and satisfies the following property. 
For any  objects $U,V$ and $W$ and any morphisms $F: U \to V$ and $G: V 
\to W,$  we 
have $*F:V \to U,$ and  the relations  
$*\textrm{Id}_U=\textrm{Id}_U$ and $*(*F) 
= F$ hold. The functor $*$ is  
called the involution or the dagger.  The morphism  $*F$ is denoted by  
$F^*,$ 
and it is called  the adjoint of  $F.$ See Burgin \cite{Burgin} or 
Brinkmann, Puppe \cite{Brink}.

Let us give some  examples of categories which are additive and dagger. 

\bigskip

\noindent {\bf Example 1:}
\begin{itemize}
\item[1)] The category of finite dimensional inner product spaces over 
$\mathbb{R}$ or $\mathbb{C}$ 
and linear maps is an example of an 
additive and a dagger category. The addition (product) of objects
is given by the  orthogonal sum and  the addition of morphism is the 
standard 
addition of linear maps.  The involution is defined as the 
adjoint of maps with 
respect to the inner products.
The existence of the adjoint to any linear map is based on the Gram-Schmid 
process which guarantees the existence of an orthonormal basis.  
The 
matrix of the adjoint of a morphism with respect to orthonormal bases in the 
domain and target spaces is given by taking the 
transpose or the transpose and the complex conjugate of the matrix of the 
original map. Such an adjoint operation on  morphisms is easily proved to be 
unambiguous.

\item[2)] The category of  
Hilbert spaces and continuous maps equipped with the addition of 
objects and maps and the involution 
given  as in item 1 is an example of an 
additive and dagger category. For the existence of the adjoints, see
Meise, Vogt \cite{Meise}. The proof is based on the Riesz representation 
theorem for Hilbert spaces.
\end{itemize}

\noindent {\bf Definition 1:} Let $\mathfrak{C}$ be an additive and dagger 
category.
We say that the {\it Hodge theory holds} for a complex $d^{\bullet}=(U^i, 
d_i: U^i \to 
U^{i+1})_{i\in \mathbb{Z}}$ in 
 $\mathfrak{C}$ or that $d^{\bullet}$   is of {\it Hodge type} if for each 
$i\in \mathbb{Z},$ we 
have
$$U^i = \textrm{Im}\ d_{i-1} \oplus \textrm{Im}\, d_i^* \oplus 
\textrm{Ker} 
\, \Delta_i$$ 
where $\Delta_i = d_i^* d_i + d_{i-1} d_{i-1}^*,$ and $d_i^*$ and 
$d_{i-1}^*$ are the adjoints
of  $d_i$ and $d_{i-1},$ respectively. We call the morphism $\Delta_i$
the $i^{th}$ {\it Laplace operator} of $d^{\bullet},$ $i\in \mathbb{Z}.$
We say that the Hodge theory holds for a subset $\mathfrak{K} \subseteq 
\mathcal{K}(\mathfrak{C})$ of complexes
in $\mathfrak{C}$ if it holds for each element $d^{\bullet} \in \mathfrak{K}.$

\bigskip

\noindent {\bf Remark 1:}
\begin{itemize}
\item[1)]
In Definition 1, we demand no compatibility of the involution 
with the additive structure. However, in the categories of pre-Hilbert and 
Hilbert $A$-modules that we will consider mostly, the relations 
$(F+G)^* = F^* + G^*$ and $(zF)^*=z^*F^*$ are satisfied for each objects $U,V,$ 
morphisms $G, F: U \to V,$ and complex number $z\in \mathbb{C}.$
\item[2)]
The existence of the Laplace operators of $d^{\bullet}$
is guaranteed by the 
definitions of the  additive and of the dagger category. 
If the dagger structure is compatible with the additive structure in the sense 
of item 1, we see that the Laplace operators are self-adjoint, i.e., 
$\Delta_i^* = \Delta_i,$ $i\in \mathbb{Z}.$
\end{itemize}

\noindent {\bf Lemma 1:} Let $d^{\bullet}=(U^i, d_i)_{i\in \mathbb{Z}}$ be a 
complex in 
the category of  Hilbert spaces  and continuous maps. If the 
images of the Laplace operators of $d^{\bullet}$ are  closed,  the Hodge 
theory 
holds for $d^{\bullet}$.

{\it Proof.} On the level of symbols, we do not distinguish 
the dependence of the inner products on 
the Hilbert spaces  and denote them by $(,).$   It is 
easy to realize 
that $\textrm{Ker}\, \Delta_i = \textrm{Ker}\, d_{i-1}^* \cap \textrm{Ker}\, 
d_i.$ Namely, the inclusion $\mbox{Ker}\, \Delta_i \supseteq \mbox{Ker}\, 
d_i 
\cap \mbox{Ker}\, d_{i-1}^*$ is immediate due to the definition of 
$\Delta_i,$ and the opposite one can be seen
as follows. For any $u \in \mbox{Ker} \, \Delta_i,$ we have $0 = (\Delta_i 
u, u) = (d_i^* d_i u  +  d_{i-1}d_{i-1}^* u,u) = 
(d_i u, d_i u) + (d_{i-1}^*u, d_{i-1}^* u).$ Since inner products are 
positive definite, we have $d_i u =0$ and 
$d_{i-1}^*u=0.$
Because  we assume the image of $\Delta_i$ to be   closed, taking the 
orthogonal complement of $\textrm{Ker} \,  \Delta_i  = \textrm{Ker}\, d_i 
\cap \textrm{Ker}\, d_{i-1}^*$,
we get  $(\textrm{Ker}\, d_{i-1}^*)^{\bot} \subseteq (\textrm{Ker}\, 
\Delta_i)^{\bot} 
= \overline{\textrm{Im}\, \Delta_i} = \textrm{Im}\, \Delta_i $ and 
$(\textrm{Ker}\, d_i)^{\bot} \subseteq   (\textrm{Ker}\, \Delta_i)^{\bot}
= \overline{\textrm{Im}\, \Delta_i} = \textrm{Im}\, \Delta_i.$ 
Summing-up, 
$$(\textrm{Ker}\, d_{i-1}^*)^{\bot} + 
(\textrm{Ker}\, d_i)^{\bot} \subseteq \textrm{Im}\, \Delta_i.$$

Further, it is immediate to see that $\textrm{Im}\, d_{i-1} \subseteq 
(\textrm{Ker}\, d_{i-1}^*)^{\bot}$ and $\textrm{Im}\, d_i^* \subseteq 
(\textrm{Ker}\, d_i)^{\bot}.$ Indeed, for any $u \in \textrm{Im}\, d_{i-1}$
there exists an element $u'\in U^{i-1}$ such that $u=d_{i-1}u'.$ For  each $v 
\in 
\textrm{Ker}\, d_{i-1}^*,$ we have $(u,v)=(d_{i-1}u',v)=(u',d_{i-1}^*v) = 0.$
Thus, the inclusion follows. The other inclusion can be seen similarly.
Using the result of the previous paragraph, 
we obtain 
\begin{equation}
\textrm{Im}\, d_{i-1} + \textrm{Im}\, d_i^* \subseteq 
(\textrm{Ker}\, d_{i-1}^*)^{\bot} + (\textrm{Ker}\, d_i)^{\bot} \subseteq 
\textrm{Im}\, \Delta_i.
\end{equation}
We prove that the sum $\textrm{Im}\, d_{i-1} + 
\textrm{Im}\, d_i^*$ is 
direct. For it,  we take $u = d_{i-1}u'$ and $v = d_i^*v'$ for $u'\in 
U^{i-1}$ and $v' \in U^{i+1},$ and compute
$(u,v)= (d_{i-1}u', d_i^*v') = (d_i d_{i-1}u',v') = 0$ which holds since 
$d^{\bullet}$ 
is a  complex. 
Therefore, we have $\textrm{Im}\, d_i^* \oplus \textrm{Im}\, 
d_{i-1}\subseteq \textrm{Im}\, \Delta_i.$
The inclusion $\textrm{Im}\, \Delta_i \subseteq 
\textrm{Im}\, 
d_{i-1} \oplus \textrm{Im}\, d_{i}^*$ is immediate. Thus, we 
conclude that  $\textrm{Im}\, \Delta_i  = \textrm{Im}\, 
d_{i-1} \oplus \textrm{Im}\, d_{i}^*.$  

Since for each $i\in \mathbb{Z},$ $\Delta_i$ is self-adjoint and its  
image 
is closed, we have    $U^i = \textrm{Im}\, \Delta_i \oplus 
\textrm{Ker} \, \Delta_i.$  
Substituting the equation for $\textrm{Im}\, \Delta_i$ found at end of the 
previous paragraph,  we get 
$U^i =  \textrm{Im}\, d_i^* \oplus \textrm{Im}\, d_{i-1}  \oplus 
\textrm{Ker} \, \Delta_i$ proving  that  the Hodge theory holds for
$d^{\bullet}.$
\hfill $\Box$

\bigskip

\noindent {\bf Remark 2:} 
 
By Lemma 1,  the Hodge theory holds for any 
complex in the category $\mathfrak{C}=V_{\textrm{fin}}$
of finite dimensional inner product spaces over real or complex numbers and 
linear maps since any 
linear subspace of a finite dimensional vector space  is closed. However, it 
is possible to prove 
that the Hodge theory holds for $\mathfrak{K} = 
\mathcal{K}(\mathfrak{C})$ in a simpler 
way than in the general case of Hilbert spaces. The relation 
$\textrm{Ker}\, 
\Delta_i = 
\textrm{Ker}\, d_i \cap \textrm{Ker} \, d_{i-1}^{*}$ is proved in the same 
way as in  the proof of Lemma 1.  
Since for any $A, B \subseteq U^i,$ the equation $(A\cap B)^{\bot} = 
A^{\bot} + B^{\bot}$ holds, we have $(\textrm{Ker}\, d_i \cap 
\textrm{Ker}\, d_{i-1}^*)^{\bot} = (\textrm{Ker}\, d_i)^{\bot} + 
(\textrm{Ker}\, d_{i-1}^*)^{\bot}.$
Due to the finite dimension, we can write
 $(\textrm{Ker}\, d_i)^{\bot} = 
\textrm{Im}\,d_{i}^*$ and $(\textrm{Ker}\, d_{i-1}^*)^{\bot} = 
\textrm{Im}\, d_{i-1},$ and thus $(\textrm{Ker}\, \Delta_i)^{\bot} = 
(\textrm{Ker}\, d_{i-1} \cap \textrm{Ker}\, 
d_i^*)^{\bot} =
\mbox{Im}\, d_{i-1} + \textrm{Im}\, d_i^*.$ The sum is   direct
as follows from $0=(d_i d_{i-1} u,v) = (d_{i-1}u, d_i^*v),$ $u\in U^{i-1},$ 
$v\in U^{i+1}.$
Substituting $(\textrm{Ker}\, \Delta_i)^{\bot} = \textrm{Im}\, d_{i-1} 
\oplus \textrm{Im}\,d_i^*$ into  $U_i = \mbox{Ker} \, 
\Delta_i \oplus (\mbox{Ker} \, \Delta_i)^{\bot},$ we get 
$U^i = \textrm{Im}\,{d}_i^* \oplus \textrm{Im}\, d_{i-1} \oplus \textrm{Ker}\, 
\Delta_i.$ 
  (Let is notice that in Lemma 1, we  proved that  the images of $d_i$ and 
$d_{i-1}^*$ are closed.)

\bigskip
 Next we define the pre-Hilbert and Hilbert modules over  $C^*$-algebras.  Our 
reference for  $C^*$-algebras is Dixmier \cite{Dixmier}.

\bigskip

\noindent {\bf Definition 2:}
For a $C^*$-algebra $A$, a {\it pre-Hilbert} $A${\it-module} is a complex 
vector space $U,$ which is a right $A$-module (the action is denoted by a dot) 
and which is moreover 
equipped  with a map
$(,):U \times U \to A$ such that for each $z\in \mathbb{C},$ $a\in A$ and 
$u,v,w\in U$ the following relations hold
\begin{itemize}
\item[1)] $(u, zv+w) = z(u,v) + (u,w)$
\item[1)] $( u, v \cdot a) = (u,v)a$ 
\item[2)] $(u,v)= (v,u)^*$
\item[3)] $(u,u) \geq 0$ and
\item[4)] $(u,u) = 0$ implies $u=0$
\end{itemize}
where $z^*$ denotes the complex conjugate of the element $z\in 
\mathbb{C}.$
A pre-Hilbert $A$-module $(U,(,))$ is called  a {\it Hilbert $A$-module} if $U$ 
is a Banach space with respect to 
the norm $U\ni u \mapsto |u| =\sqrt{|(u,u)|_A}\in [0,+\infty).$ The map $(,): 
U\times U \to A$ is called the  {\it $A$-product}. 
\bigskip

Note that if $A$ is the  algebra of complex numbers, Definition 2 coincides 
with the one of a pre-Hilbert and of a Hilbert  space, respectively.

{\it Morphisms} of pre-Hilbert $A$-modules $(U,(,))$ and $(V,(,)_V)$ 
are assumed to be  continuous, $A$\--li\-ne\-ar and adjointable  
maps. Recall that a map $L:U \to V$ is called $A$-linear if the  equivariance 
condition
$L(u)\cdot a = L(u \cdot a)$ holds for any  $a\in A$ and $u\in U.$  
An adjoint  $L^*: V \to U$ of a pre-Hilbert $A$-module morphism $L: U \to V$
is a map which satisfies $(Lu, v)_V = (u,L^*v)_U$ for any $u\in U$ and $v\in 
V.$ It 
 is known that the adjoint need not exist in general, and that if it exists,
it is unique and a pre-Hilbert $A$-module homomorphism, i.e., continuous and 
$A$-linear.
Morphisms of Hilbert $A$-modules have to be morphisms 
of these modules considered as pre-Hilbert $A$-modules. 
The category the objects of which are pre-Hilbert $A$-modules and the morphisms 
of which are continuous, $A$-linear and adjointable maps will be denoted by
 $PH_A^*.$  The category $H_A^*$ of Hilbert $A$-modules   is defined to be the 
full subcategory of $PH_A^*$ the object of which are Hilbert $A$-modules.
If we drop the condition on the adjointability of morphisms, we denote the 
resulting categories by $PH_A$ and $H_A.$
By an {\it isomorphism} $F: U \to V$ in $PH_A^*$ or $H_A^*,$ 
we mean  a morphism which is right and left invertible by a morphism in 
$PH_A^*$ or $H_A^*,$ respectively.
In particular, we demand an isomorphism in these categories  neither to 
preserve the appropriate
$A$-products nor the induced norms.  

Submodules  of a (pre-)Hilbert $A$-module have to be (pre-)Hilbert $A$-modules 
with respect to the restrictions both of the 
algebraic and of the norm structure. In particular, they are 
closed in the super-module. 
Further, if $U$ is
a submodule of the (pre-)Hilbert $A$-module  $V,$  we can construct the space 
$U^{\bot} = \{ v\in V, \, (v,u) = 0 \textrm{ for all } u\in U\}$ which is a 
(pre-)Hilbert $A$-module.
Further, $U$ is called orthogonally complemented in $V$ if
$V = U\oplus U^{\bot}.$ There are Hilbert $A$-submodules which are not 
orthogonally complemented. (See Lance \cite{Lance}.)
For the convenience of the reader, we give several examples of  Hilbert 
$A$-modules 
and an example of a pre-Hilbert $A$-module.
For further examples, see   Solovyov, Troitsky 
\cite{ST},  Manuilov, Troitsky \cite{MT}, Lance \cite{Lance}, and 
Wegge-Olsen \cite{Wegge}.

\bigskip

\noindent{\bf Example 2:} 
\begin{itemize}
\item[1)] Let $H$ be a Hilbert space with  the inner product denoted by 
$(,)_H.$ The  
 action of the $C^*$-algebra $A=B(H)$ of bounded linear operators on $H$
is by evaluation on the adjoint, i.e., $h\cdot a =a^*(h)$ for any $a\in B(H)$ 
and $h\in H.$ The $B(H)$-product  is defined by 
$(u,v) = u \otimes v^*,$ where
$(u\otimes v^*)w = (v,w)_Hu$ for $u,v,w \in H.$ 
In this case, the product takes values in  the $C^*$-algebra $K(H)$ of compact 
operators on $H.$ In fact, the $A$-product maps into the 
algebra of finite rank operators.

\item[2)] For a locally compact topological space $X,$ consider the 
$C^*$-algebra 
$A=\mathcal{C}_0(X)$ of continuous functions vanishing at infinity with the 
product given by the point-wise multiplication, 
with the  complex conjugation as the involution, and with the classical 
supremum 
norm 
$|\,|_A:\mathcal{C}_0(X) \to [0,+\infty)$
$$|f|_A = \textrm{sup}\{|f(x)|,\, x\in X\}$$ where $f\in A.$ For  $U$, we take 
the $C^*$-algebra $\mathcal{C}_0(X)$ itself with the module 
structure given by the point-wise multiplication, i.e., $(f\cdot g)(x) = 
f(x)g(x),$ $f\in U,$ $g\in A$ and $x\in X.$
The $A$-product is defined by $(f,g) = \overline{f}g.$  Note that this is a 
particular 
example of a Hilbert $A$-module with $U=A,$ right action
$a\cdot b = ab$ for $a\in U=A$ and $b\in A,$ and $A$-product $(a,b) = 
a^*b,$ $a,b \in U.$

\item[3)]If $U$ is a Hilbert $A$-module,  the 
 orthogonal direct sums of a finite number of copies of $U$ form a Hilbert 
$A$-module in a 
natural way.
One can also construct the space $\ell^2(U),$  i.e., the space consisting of 
sequences 
$(a_n)_{n \in \mathbb{N}}$ with $a_n \in U,$ $n\in \mathbb{N},$ for 
which the series $\sum_{i=1}^{\infty}(a_i,a_i)$  converges in $A$. The 
$A$-product is given by  $((a_n)_{n\in \mathbb{N}}, (b_n)_{n\in \mathbb{N}}) 
= \sum_{i=1}^{\infty}(a_i, b_i),$ where $(a_n)_{n \in \mathbb{N}}, (b_n)_{n\in 
\mathbb{N}} \in \ell^2(U).$
See Manui\-lov, Troitsky \cite{MT}.

\item[4)] Let $A$ be a $C^*$-algebra. For a compact manifold $M^n,$ pick a 
Riemannian metric $g$ and choose a volume element 
$|\textrm{vol}_g|\in \Gamma(M,|\bigwedge^nT^*M|).$ Then for any $A$-Hilbert 
bundle $\mathcal{E} \to M$ with fiber a Hilbert $A$-module $E$,
 one defines a {\it pre-Hilbert $A$-module}
$\Gamma(M,\mathcal{E})$ of smooth sections of $\mathcal{E} \to M$ by setting
 $(s\cdot a)_m= s_m \cdot a$ for  $a\in A,$  $s\in \Gamma(M,\mathcal{E}),$ and 
$m\in M.$
One sets 
$$(s',s)= \int_{m\in M}(s'_m, s_m)_m |\textrm{vol}_g|_m$$ where $s,s'\in 
\Gamma(M,\mathcal{E})$, $(,)_m$ denotes the $A$-product in  fiber 
$\mathcal{E}_m,$ and $m\in M.$
Taking the  completion of $\Gamma(M,\mathcal{E})$ with respect to the norm 
 associated to the $A$-product $(,)$ (Definition 2), we get the Hilbert 
$A$-module $(W^0(M,\mathcal{E}), (,)_0).$
Further Hilbert $A$-modules $(W^t(M, \mathcal{E}),(,)_t),$ $t\in \mathbb{N}_0$ 
are 
derived
from the space $\Gamma(M,\mathcal{E})$ by mimicking  the construction of  
Sobolev  
spaces  defined   for  finite rank bundles. See Wells \cite{Wells} for the 
finite rank case and  Solovyov, Troitsky \cite{ST} for the case of $A$-Hilbert 
bundles. 
\end{itemize}

Let us turn our attention to the so-called self-adjoint parametrix possessing 
morphisms in the category $\mathfrak{C}=PH_A^*.$

\bigskip

\noindent {\bf Definition 3:} A pre-Hilbert $A$-module endomorphism $F: U \to 
U$ is 
called 
{\it self-adjoint parametrix possessing} if
$F$ is self-adjoint, i.e., $F^*=F,$ and there exist a pre-Hilbert $A$-module 
homomorphism $G: U\to U$ and
a self-adjoint pre-Hilbert $A$-module homomorphism 
$P: U \to U$  such that
\begin{eqnarray*} 
1_U  &=& GF + P\\ 
1_U  &=& FG + P \\
FP &=& 0.
\end{eqnarray*}

\bigskip

\noindent {\bf Remark 3:} 
\begin{itemize}
\item[1)] The map $G$ from Definition 3 is called a {\it parametrix} or a 
{\it Green operator} and the first two equations in this definition
are called the {\it parametrix equations}.
\item[2)]
Composing the first  parametrix equation from the right with $P$ and using the 
third equation, we get that $P^2 = P.$
\item[3)] If $F:U \to U$ is a self-adjoint parametrix possessing morphism in 
$PH_A^*$, then $U=\textrm{Ker}\, F 
\oplus \textrm{Im} \, F$ (see Theorem 6 in Kr\'ysl \cite{KryslAGAG2}).  
In particular, the image of $F$ is closed.  Note that we do not assume that
$U$ is complete.
\item[4)]A morphism in $H_A^*$ is self-adjoint 
parametrix possessing if its image is closed.
Indeed, the  Mishchenko theorem (Theorem 3.2 on pp. 22 in Lance 
\cite{Lance}) enables us to write for a self-adjoint morphism $F: U \to U$ 
with closed image, the 
orthogonal decomposition  $U = \mbox{Ker}\, F \oplus \mbox{Im}\, F.$
Then we can define the projection onto $\mbox{Ker}\, F$ along $\mbox{Im}\, F.$
It is easy to see that the projection is 
self-adjoint.  Inverting $F$ on its image and defining it by zero on the kernel 
of $F$,  we get a map $G$ which satisfies the parametrix equations and it is 
continuous 
due to the open map theorem.  Thus, in $H_A^*$ a self-adjoint map $F$ 
is self-adjoint parametrix possessing if and only if its image is closed.
\end{itemize}

Let us notice that if  $d^{\bullet} = (U^i, d_i)_{i \in \mathbb{Z}}$ is a 
co-chain complex in the category $PH_A^*,$ the 
$i^{th}$ Laplace operator $\Delta_i  = d_{i-1}d_{i-1}^* + d_i^* d_i$ is  
self-adjoint, $i\in \mathbb{Z}.$

\bigskip

\noindent {\bf Definition 4:}
A co-chain complex $d^{\bullet} \in \mathcal{K}(PH_A^*)$ is called {\it 
self-adjoint parametrix possessing}
if all of its Laplace operators are self-adjoint parametrix possessing maps.

\bigskip

\noindent{\bf Remark 4:}
\begin{itemize}
\item[1)] Since 
$\Delta_{i+1}d_i = (d_{i+1}^*d_{i+1} + d_id_i^*)d_i = d_i d_i^* d_i =  
d_id_i^*d_i + d_id_{i-1}d_{i-1}^*
= d_i(d_i^*d_i + d_{i-1}d_{i-1}^*) = d_i \Delta_i,$
the Laplace operators are co-chain 
endomorphisms of $d^{\bullet}.$
Similarly, one derives that the Laplace operators are
chain endomorphisms of the chain complex $(U^i, d_i^*: U^{i+1} \to 
U^i)_{i\in 
\mathbb{Z}}$  "dual" to $d^{\bullet}.$

\item[2)] Let us assume that the Laplace operators $\Delta_i$ of a 
complex 
$d^{\bullet}$ in $PH_A^*$ satisfy  
equations $\Delta_i G_i + P_i = G_i 
\Delta_i + P_i = 1_{U_i}$ and that the identity $\Delta_iP_i=0$ holds. Notice 
that we do not suppose that the idempotent $P_i$ is self-adjoint. Still, we 
can  prove that the
Green 
operators $G_i$ satisfy
$G_{i+1}d_i = d_i G_i,$ i.e., that they are co-chain endomorphisms of 
the complex $d^{\bullet}$ we consider.  
For it, see Theorem 3, Kr\'ysl \cite{KryslAGAG1}.    

\item[3)] In the following picture,  facts from the previous two items 
are summarized in a diagrammatic way.

 $$\xymatrix@C=2cm@R=2cm@=2cm{
U_{i-1}\ar@<2pt>[r]^{d_{i-1}}  \ar@<2pt>[d]^{G_{i-1}}  & 
U_i \ar@<2pt>[l]^{d_{i-1}^*}  \ar@<2pt>[r]^{d_i}  \ar@<2pt>[d]^{G_{i}}  & 
U_{i+1 } \ar@<2pt>[d]^{G_{i+1}} \ar@<2pt>[l]^{d_{i}^*}\\
U_{i-1} \ar@<2pt>[r]^{d_{i-1}} \ar@<2pt>[u]^{\Delta_{i-1}} &       U_i  
\ar@<2pt>[u]^{\Delta_i}  \ar@<2pt>[r]^{d_i} \ar@<2pt>[l]^{d_{i-1}^*} &  
U_{i+1} \ar@<2pt>[l]^{d_{i}^*} \ar@<2pt>[u]^{\Delta_i}
}$$

\end{itemize}

Let us consider the  cohomology groups $H^i(d^{\bullet}) = 
\textrm{Ker}\, 
d_i / \textrm{Im}\, d_{i-1}$ of a complex $d^{\bullet} \in 
\mathcal{K}(PH_A^*),$  $i \in 
\mathbb{Z}.$
If  $\textrm{Im} \, d_{i-1}$ is orthogonally complementable
in $\textrm{Ker} \, d_i,$ then one can define an $A$-product in 
$H^i(d^{\bullet})$ by $([u],[v])_{H^i(d^{\bullet})} =  (p_iu,p_iv),$ where 
$u,v\in U^i$ and
$p_i$ is the projection 
along $\textrm{Im} \, d_{i-1}$ 
onto the orthogonal complement  $(\textrm{Im}\, d_{i-1})^{\bot}$ in 
$\textrm{Ker}\,d_i.$ Let us call this $A$-product the {\it canonical 
quotient product}. For information on $A$-products on  quotients in 
$PH_A^*,$ see  \cite{KryslAGAG2}.  

\bigskip
In the next theorem, we collect results on self-adjoint parametrix complexes 
from \cite{KryslAGAG2}.
 
\bigskip

\noindent {\bf Theorem 2:} Let $A$ be a $C^*$-algebra. If
 $d^{\bullet}=(U^i,d_i)_{i \in \mathbb{Z}} \in \mathcal{K}(PH_A^*)$  is 
self-adjoint parametrix possessing complex,
then  for any $i\in \mathbb{Z},$
\begin{itemize}
\item[1)] $U^i = \textrm{Ker} \, \Delta_i \oplus \textrm{Im}\, d_i^* \oplus 
\textrm{Im} \, d_{i-1},$ i.e., $d^{\bullet}$ is a Hodge type complex
\item[2)] $\textrm{Ker}\, d_{i} = \textrm{Ker} \, \Delta_i \oplus 
\textrm{Im} 
\, 
d_{i-1}$
\item[3)] $\textrm{Ker}\, d_{i}^* = \textrm{Ker} \, \Delta_{i+1} \oplus 
\textrm{Im} 
\, d_{i+1}^*$
\item[4)] $\textrm{Im}\, \Delta_i  = \textrm{Im} \, d_i^* \oplus \textrm{Im} 
\,  
d_{i-1}$
\item[5)] $H^i(d^{\bullet})$ is a pre-Hilbert $A$-module with respect to the 
canonical quotient product $(,)_{H^i(d^{\bullet})}$
\item[6)] The spaces $\textrm{Ker} \, \Delta_i$ and  $H^i(d^{\bullet})$ 
are isomorphic as pre-Hilbert $A$-mod\-u\-les.
Moreover, if  $d^{\bullet}$ is self-adjoint parametrix possessing com\-plex in 
$\mathcal{K}(H_A^*),$ then $H^i(d^{\bullet})$ is an $A$-Hilbert module
and $\textrm{Ker} \, \Delta_i \simeq H^i(d^{\bullet})$ are isomorphic as 
$A$-Hilbert modules.
\end{itemize}

{\it Proof.} See Theorem 11 in Kr\'ysl \cite{KryslAGAG2} for item 1; 
Theorem 13 in \cite{KryslAGAG2} for items 2 and 3; Remark 12 (1) in 
\cite{KryslAGAG2} for item 4; and 
Corollary 14 in \cite{KryslAGAG2} for items 5 and 6.  
\hfill $\Box$

\bigskip
Next we prove that in the category $\mathfrak{C}=H_A^*,$ the property of a 
complex to be self-adjoint parametrix possessing {\bf characterizes} the 
complexes of Hodge type.

\bigskip
\noindent {\bf Theorem 3:}
If the Hodge theory holds for a complex $d^{\bullet} \in \mathcal{K}(H_A^*),$ 
then $d^{\bullet}$ is self-adjoint parametrix possessing.

{\it Proof.}
 Because the Hodge theory holds for $d^{\bullet},$ we have the decomposition of 
$U^i$ into
	Hilbert $A$-modules
  \begin{align*}
  U^i = \textrm{Ker}\,\Delta_i \oplus \textrm{Im}\, d_{i-1} \oplus 
\textrm{Im} \, d_i^*
\end{align*} $i\in \mathbb{Z}.$ In particular, the ranges of $d_{i-1}$ and 
$d_i^*$ are 
closed topological vector spaces.
It is easy to verify that 
\begin{align*}
\textrm{Ker}\, d_i^*d_i = \textrm{Ker}\, d_i && \textrm{Ker}\,
d_{i-1}^* = \textrm{Ker} \, d_{i-1}d_{i-1}^*.
\end{align*}
  For $i \in \mathbb{Z}$ and $u\in U^i,$ we have 
$$(\Delta_i u, \Delta_i u) = (d_i^*d_i u,d_i^* d_iu) + (d_{i-1}d_{i-1}^* u, d_{i-1} d_{i-1}^* u)$$ 
since $(d_i^*d_iu, d_{i-1}d_{i-1}^*u) = (d_i u, d_i d_{i-1}d_{i-1}^* u) = 0$ for any
$u\in U^i.$ 
Due to the definition of the Laplace operator and the positive definiteness of the $A$-Hilbert product,
we have $\textrm{Ker}\, \Delta_i = \textrm{Ker}\, d_i \cap \textrm{Ker}d_{i-1}^*.$
For $u \in (\textrm{Ker}\,\Delta_i)^{\bot} = (\textrm{Ker}\, d_i)^{\bot} + 
(\textrm{Ker} \, d_{i-1}^*)^{\bot},$
there exist $u_1 \in (\textrm{Ker}\, d_i)^{\bot} = \textrm{Im}\, d_i^* 
$ and $u_2 \in (\textrm{Ker} \, d_{i-1}^*)^{\bot}= \textrm{Im}\, 
d_{i-1}$ such that $u=u_1+u_2.$
Consequently, $(\Delta_i u, \Delta_i u) = $
\begin{eqnarray*}
                         &=& (d_i^*d_i(u_1+u_2),d_i^*d_i (u_1+u_2)) +  (d_{i-1}d_{i-1}^*(u_1+u_2), d_{i-1}d_{i-1}^*(u_1+u_2)) \\
                         &=& (d_i^*d_iu_1,d_i^*d_iu_1) + (d_i^*d_iu_2,d_i^*d_i  
u_2) + (d_i^*d_iu_1,d_i^*d_iu_2) + (d_i^*d_i u_2,d_i^*d_i u_1) \\
                        & & + (d_{i-1}d_{i-1}^*u_1, d_{i-1}d_{i-1}^*u_1) + 
(d_{i-1}d_{i-1}^*u_2, d_{i-1}d_{i-1}^*u_2)  \\
& & +(d_{i-1}d_{i-1}^*u_1, d_{i-1}d_{i-1}^*u_2) + (d_{i-1}d_{i-1}^*u_2, 
d_{i-1}d_{i-1}^*u_1)\\
&=& (d_i^*d_iu_1,d_i^*d_iu_1) + (d_{i-1}d_{i-1}^*u_2, 
d_{i-1}d_{i-1}^*u_2)
\end{eqnarray*}
since $(d_i^*d_iu_2,d_i^*d_i 
u_2)= (d_i^*d_iu_1,d_i^*d_iu_2)=0$ due to $u_2 \in \textrm{Im}\, 
d_{i-1},$ and \\ $(d_{i-1}d_{i-1}^*u_1, 
d_{i-1}d_{i-1}^*u_1)=(d_{i-1}d_{i-1}^*u_1, d_{i-1}d_{i-1}^*u_2)=0$ due to 
$u_1 \in \textrm{Im}\, d_i^*.$
Since both summands on the right-hand side of $$(\Delta_iu,\Delta_iu) = 
(d_i^*d_iu_1,d_i^*d_iu_1) + (d_{i-1}d_{i-1}^*u_2, 
d_{i-1}d_{i-1}^*u_2)$$ are 
non-negative, we obtain
$(\Delta_i u, \Delta_i u)\geq  (d_i^*d_iu_1, d_i^*d_iu_1)$ and $(\Delta_i u, 
\Delta_i u) \geq (d_{i-1}d_{i-1}^*u_2, d_{i-1}d_{i-1}^*u_2).$
Consequently
\begin{eqnarray}
|\Delta_i u| &\geq& |d_i^*d_i u_1|\\
|\Delta_i u| &\geq& |d_{i-1}d_{i-1}^* u_2|
\end{eqnarray}
(See  paragraph  1.6.9 on pp. 18 in  Dixmier \cite{Dixmier}.)
Notice that $d_i^*d_i$ and $d_{i-1}^*d_{i-1}$ is injective on 
$(\textrm{Ker}\, d_i^*d_i)^{\bot} = (\textrm{Ker}\, d_i)^{\bot}$ and 
$(\textrm{Ker}\, d_{i-1}d_{i-1}^*)^{\bot}=\\=(\textrm{Ker}\, 
d_{i-1}^*)^{\bot},$  respectively, and zero on the complements of these spaces. 
	Due to 
an equivalent 
characterization of closed image maps on  Banach spaces,
there are positive real numbers
$\alpha, \beta$ such that $|d_i^*d_i u_1| \geq \alpha|u_1|$ and $|d_{i-1}d_{i-1}^*u_2| \geq \beta |u_2|$ (see, e.g., Abramovich, Aliprantis \cite{Abra}).
Substituting these inequalities into (2) and (3) and adding the resulting 
inequalities,  we see that
$2|\Delta_i u| \geq \alpha|u_1| + \beta|u_2|.$
Thus $|\Delta_i u| \geq \frac{1}{2}\textrm{min}\{\alpha,\beta\}
(|u_1| +|u_2|)  \geq \frac{1}{2}\textrm{min}\{\alpha,\beta\} |u_1 + 
u_2|=\frac{1}{2}\textrm{min}\{\alpha,\beta\} |u|$ by the triangle identity.
Due to the characterization of closed image maps again, we get that the image of 
$\Delta_i$ is closed.
This implies that $d^{\bullet}$ is self-adjoint parametrix possessing using 
Remark 3 item 4.

\hfill $\Box$

\bigskip

\noindent {\bf Remark 5:}
\begin{itemize}
\item[1)]
Let $\mathfrak{C}= H_{\mathbb{C}}^*$ be the  category of Hilbert spaces and 
continuous maps.  Let us consider 
such complexes in $\mathfrak{C}$ whose differentials are Fredholm maps.
Especially their images and the ones of their adjoint maps are closed.
It is easy  to prove that 
the Laplacians are Fredholm as well. Namely, due to the fact that the images of 
$d_i$ and $d_i^*,$ $i\in \mathbb{Z},$ are closed we can use  the 
inequalities in the proof of Theorem 3 and conclude that the Laplacian has 
closed image by the characterization of closed image maps as used above. By 
Lemma 1, the complex is necessarily of Hodge type  and we have
the decomposition $\textrm{Im}\, \Delta_i = \textrm{Im}\, d_{i-1} \oplus 
\textrm{Im}\, d_i^*$ (e.g., by Theorem 2). In particular, the cokernel of the 
Laplacian is a finite dimensional space. The kernel of $\Delta_i$ is finite 
dimensional due to the finite dimensionality of $\textrm{Ker}\, d_i$ and due to
$\textrm{Ker}\, \Delta_i = \textrm{Ker} \, d_i \cap \textrm{Ker}\, d_{i-1}^*$ 
which follows from the definition of the Laplacian operator.
\item[2)] From Theorem 2,  Theorem 3 and Remark 3 item 4, we get that a complex 
in $H_A^*$ is of Hodge type if and only if the 
images of its Laplace operators are closed if and only if it is self-adjoint 
parametrix possessing.
\end{itemize}

\bigskip
 
\noindent {\bf Example 3:}
\begin{itemize}
\item[1)]For a compact manifold  $M$  of positive dimension,  
let
us consider the Sobolev spaces $W^{k,l}(M)$ for $k,l$ non-negative integers. 
For $l=2,$ these space are complex Hilbert spaces.
Due to the Rellich-Kondrachov embedding theorem 
and the fact that the
dimension of $W^{k,2}(M)$ is infinite, 
the  canonical embedding $i:W^{k,2}(M) \hookrightarrow 
W^{l,2}(M)$ has a non-closed image for 
$k>l.$ 
We take  $d^{\bullet}=$ 
$$\xymatrix{
 0  \ar[r] &   W^{k,2}(M) \ar[r]^{i} &  W^{l,2}(M) 
\ar[r]&    0}.$$ Labeling the first element in the complex by zero, 
the second cohomology $H^2(d^{\bullet}) = \textrm{Ker} 
\, 0/\textrm{Im} \, i =
W^{l,2}(M)/i(W^{k,2}(M))$ is non-Hausdorff in the quotient topology.
The complex is not self-adjoint parametrix possessing due to 
Theorem 2 item 5. Consequently, it is not of Hodge type (Theorem 3).
 
\item[2)] This example shows a simpler construction of a complex 
in $\mathcal{K}(H_{\mathbb{C}}^*)$ which is not of Hodge type. Without any 
reference to a manifold, we can define  mapping $i: 
\ell^2(\mathbb{Z})\to 
\ell^2(\mathbb{Z})$ by setting
$i(e_n) = e_n /n,$ where $(e_n)_{n=1}^{+\infty}$ denotes the canonical 
orthonormal 
system of $\ell^2(\mathbb{Z}).$
It is easy to check that  $i$ is continuous. Further, the set  
$i(\ell^2(\mathbb{Z}))$ 
is not closed.
For it, the sequence $(1, 1/2, 1/3,\ldots ) \in \ell^2(\mathbb{Z})$ is 
not in the image.
Indeed, the preimage of this element had to be the sequence $(1,1,1,\ldots)$ 
which is not in $\ell^2(\mathbb{Z}).$ On the other hand,
$(1, 1/2, 1/3,\ldots )$ lies in the closure of $i(\ell^2(\mathbb{Z}))$ since it 
is the limit of 
the sequence $i((1,0\ldots)), i((1,1,0, \ldots)), i((1,1,1,0\ldots)), \ldots .$ 
 The complex $0 \rightarrow \ell^2(\mathbb{Z}) \xrightarrow{i} 
\ell^2(\mathbb{Z}) \rightarrow  0$ is not of Hodge type and it is not 
self-adjoint parametrix 
possessing by similar reasons as given in the example above.
\end{itemize}

\section{$C^*$-Fredholm operators  over 
$C^*$-algebras \\ of compact 
operators}

In this section, we focus on complexes over $C^*$-algebras of compact 
operators, and study $C^*$-Fredholm maps acting between Hilbert modules over 
such algebras.  For the convenience of the reader, let us  recall some 
necessary notions.  

\bigskip

\noindent {\bf Definition 5:} 
Let $(U,(,)_U)$ and $(V,(,)_V)$ be Hilbert $A$-modules.
\begin{itemize}
\item[1)] For any $u\in U$ and $v\in V,$  the operator $F_{u,v}: U \to V$  
defined by $U 
\ni u' \mapsto F_{u,v}(u') = v\cdot(u,u')$ is called an {\it elementary 
operator}.
A morphism  $F: U \to V$ in $H_A^*$ is called of {\it $A$-finite rank} if it 
can 
be written as a finite sum of $A$-linear combinations
of the elementary operators. 
\item[2)] The set $K_A(U,V)$ of $A$-{\it compact} operators on $U$ is defined 
to be the closure of  the vector space  of the $A$-finite rank morphisms in the 
operator norm in $\textrm{Hom}_{H_A^*}(U,V),$ induced by the norms $|\,|_U$ and 
$|\,|_V.$ 

\item[3)] We call  $F \in \textrm{Hom}_{H_A^*}(U,V)$ $A$-{\it Fredholm} if 
there exist Hilbert $A$-module homomorphisms $G_V:V \to U$ and $G_U:U \to V$ 
and $A$-compact homomorphisms $P_U:U \to U$ and $P_V:V \to V$ such that 
$$G_UF=1_U + P_U $$ 
$$FG_V = 1_V + P_V$$ i.e.,
if  $F$ is left and right invertible modulo $A$-compact operators.
\end{itemize}

\bigskip

\noindent {\bf Remark 5:}
\begin{itemize}
\item[1)] {\it Equivalent definition of $A$-compact operators.} The
 $A$-finite rank operators  are easily seen to be 
adjointable. 
Suppose for a moment that we define  the "$A$-compact" operators as such 
morphisms in the category $H_A$ of Hilbert $A$-modules and continuous
$A$-module homomorphisms  that lie 
in the operator norm closure in
$\textrm{Hom}_{H_A}(U,V) \supseteq \textrm{Hom}_{H_A^*}(U,V)$ of the $A$-finite 
rank operators. 
One can prove that these operators  are  adjointable and that their set
coincides with class of 
the $A$-compact operators defined above (Definition 5).
For it see, e.g., Corollary 15.2.10 in Wegge-Olsen \cite{Wegge}.

\item[2)] {\it $A$-compact vs. compact.} It is well known that in general, the 
notion of an $A$-compact operator does not coincide 
with the notion of a compact operator in a Banach space. Indeed, let us 
consider an infinite dimensional {\it unital}
$C^*$-algebra $A$ $(1\in A),$ and take $U=A$ with the right action given by the 
multiplication in $A$ and the $A$-product $(a,b) = a^*b,$ $a, b \in A.$
Then the identity  $1_U:U \to U$ is $A$-compact since it is equal to 
$F_{1,1}.$ But it is not a compact operator in the classical  sense since  $U$ 
is infinite 
dimensional. 

\end{itemize}

\bigskip
\noindent {\bf Example 4:}

\begin{itemize}
\item[1)]
{\it $A$-Fredholm operator with non-closed image.} Let us consider the space 
$X=[0,1]\subseteq \mathbb{R},$ the $C^*$-algebra $A=\mathcal{C}([0,1])$ and the 
tautological Hilbert $A$-module $U=A=\mathcal{C}([0,1])$ 
(second paragraph of  Example 2). We give a simple proof of the
fact  that 
there exists  an endomorphism on $U$ which is $A$-Fredholm but the image of 
which is 
not  closed.  Let us take an arbitrary map $T \in \textrm{End}_{H_A^*}(U).$ 
Writing $f= 1 
\cdot f,$ we have $T(1 \cdot f) = T(1)\cdot f = T(1)f.$
 Thus, $T$ can be written as  the elementary operator $F_{f_0,1}$ where $f_0 = 
T(1).$ Since $T$ is arbitrary, $K_A(U,U)=\textrm{End}_{H_A^*}(U).$ 
Consequently, any endomorphism $T \in \textrm{End}_{H_A^*}(U)$
is $A$-Fredholm since $T1_U=1_UT=1_U+(T-1_U)$  
Let us consider   operator $Ff = xf,$ $f\in U.$ This operator satisfies 
$F = F^*,$ and it is clearly a morphism of the Hilbert $A$-module $U.$  
It is immediate to realize that $\textrm{Ker}\, F = 0.$
Suppose  that the image of $F=F^*$ is closed. Using 
Theorem 3.2 in  Lance \cite{Lance},
we  obtain $\mathcal{C}([0,1]) = \textrm{Im}\, F^* \oplus \textrm{Ker} \, F = 
\textrm{Im} \, F.$ Since the constant function $1 \notin \textrm{Im}\, F,$ we 
get a contradiction.
Therefore  $\textrm{Im}\, F$ is not closed although  $F$ is an $A$-Fredholm 
operator  
as shown  above.
Let us recall that the image of a Fredholm operator on a Banach space, in the 
classical sense, is closed.

\item[2)] {\it Hilbert space over its compact operators.}
Let us note that if $U = H$ and $A=K(H)$ with the action and the 
$A$-product as in Example 2 item 1, we have $F_{u,v}=(v,u)$ 
for any $u,v \in H.$ Especially, $F_{u,v}$ are  rank one operators. Thus, their
 finite $A$-linear combinations are finite rank operators on $H$ and their 
closure is  $K(H)$ itself, i.e., $K_{K(H)}(H)=K(H).$ 
\end{itemize}

\bigskip

\noindent {\bf Remark 6:} Let us remark that the definition of an $A$-Fredholm 
operator 
on pp. 841  in Mish\-chenko, Fo\-menko \cite{FM} is different from  the 
definition of an $A$-Fredholm 
operator given in item 3 of  Definition 5 of our paper. However, an 
$A$-Fredholm 
operator in the sense of Fomenko and Mishchenko is  necessarily invertible 
modulo an $A$-compact operator 
(see Theorem 2.4 in Fomenko, Mishchenko \cite{FM}), i.e., it is an $A$-Fredholm 
operator in our sense.

\bigskip

\noindent {\bf Definition 6:} A $C^*$-algebra is called a  $C^*${\it -algebra 
of compact 
operators} if it is a  $C^*$-subalgebra of the $C^*$-algebra of 
  compact operators $K(H)$ on a   Hilbert space $H.$ 

 \bigskip

If $A$ is a $C^*$-algebra of compact operators, an analogue of an  orthonormal 
system in a Hilbert space is introduced for the case of Hilbert $A$-modules  
in the paper of Baki\'c, Gulja\v{s} \cite{Guljas}. For a fixed
Hilbert $A$-module, the cardinality of any of its orthonormal systems does not 
depend on the  choice of such a system. 
We denote the cardinality of an orthonormal system of a Hilbert $A$-module $U$ 
over  a $C^*$-algebra $A$ of compact operators   by $\textrm{dim}_AU.$ Let us 
note that in particular, an orthonormal system forms a set of 
generators of the  module as follows from the definition in \cite{Guljas}.

\bigskip

\noindent {\bf Theorem 4:} Let $A$ be a $C^*$-algebra of compact operators, $U$ 
be a 
Hilbert $A$-module, and  $F \in \textrm{End}_{H_A^*}(U).$ 
Then $F$ is $A$-Fredholm, if and only if its image is closed and  
$\textrm{dim}_A \textrm{Ker}\, F$ and $\textrm{dim}_A(\textrm{Im} \, 
F)^{\perp}$ 
are finite.

{\it Proof.}  Baki\'c,  Gulja\v{s} \cite{Guljas}, pp. 268.
\hfill $\Box$

\bigskip

\noindent {\bf Corollary 5:} Let $A$ be a $C^*$-algebra of compact operators, 
$U$ and
$V$ be Hilbert $A$-modules, and $F \in \textrm{Hom}_{H_A^*}(U,V).$ Then
$F$ is an $A$-Fredholm operator, if and only if its image is closed and
$\textrm{dim}_A \textrm{Ker} \, F$ and $\textrm{dim}_A(\textrm{Im}\, 
F)^{\perp}$ are finite.

{\it Proof.} Let $F: U \to V$ be an $A$-Fredholm operator and
$G_U, P_U$ and $G_V, P_V$ be the corresponding left and right 
inverses and projections, respectively, i.e., 
$G_UF = 1_U + P_U$ and $FG_V = 1_V + P_V.$

Let us consider the element $\mathfrak{F}= \bigl(\begin{smallmatrix}
0 & F^* \\ F& 0
\end{smallmatrix} \bigr) \in \textrm{End}_{H_A^*}(U\oplus V).$
For this element, we can write  
$$\left(\begin{matrix} 0 &G_U \\ G_V^* & 0\end{matrix} \right) 
\left(\begin{matrix}
0 & F^* \\ F& 0
\end{matrix} \right)  = 
\left(\begin{matrix}
1_U + P_U & 0 \\ 0 & 1_V + P_V \end{matrix}\right)= \left( \begin{matrix} 1_U & 
0 \\ 0 & 1_V \end{matrix}\right) +\left( \begin{matrix} P_U & 0\\0 & P_V  
\end{matrix} \right)  $$
Since  the last written matrix is an $A$-compact operator in 
$\textrm{End}_{H_A^*}(U\oplus V),$
$\mathfrak{F}$ is left invertible modulo an $A$-compact operator on $U \oplus 
V.$ 
The right invertibility is proved in a similar way. Summing-up, $\mathfrak{F}$ 
is 
$A$-Fredholm.
According to Theorem 4, $\mathfrak{F}$ has closed image. This implies that
$F$ has closed image as well due to the orthogonality of the modules $U$ and 
$V$ in  $U \oplus V.$ Let us denote the orthogonal projections of $U\oplus V$
onto $U$ and $V$ by $\textrm{proj}_U$ and $\textrm{proj}_V,$  respectively.
Due to Theorem 4, $\textrm{dim}_A(\textrm{Ker}\, 
\mathfrak{F})$ and $\textrm{dim}_A(\textrm{Im}\, 
\mathfrak{F})^{\bot}$ are finite.
Since $\textrm{Ker}\, F= \textrm{proj}_U (\textrm{Ker}\, \mathfrak{F})$
and $(\textrm{Im}\, F)^{\bot} = \textrm{Ker} \, F^* = \textrm{proj}_V( 
\textrm{Ker}\,\mathfrak{F}^*) = \textrm{proj}_V (\textrm{Im} \, 
\mathfrak{F})^{\bot},$  the finiteness of 
$\textrm{dim}_A \textrm{Ker} \, F$ and $\textrm{dim}_A(\textrm{Im}\, 
F)^{\perp}$ follows from the finiteness of $\textrm{dim}_A(\textrm{Ker}\, 
\mathfrak{F})$ and  $\textrm{dim}_A(\textrm{Im}\, 
\mathfrak{F})^{\bot}.$

On the other hand, if $\textrm{dim}_A \, \textrm{Im}\, F$ and $\textrm{dim}_A 
(\textrm{Im}\, F)^{\bot}$  are
finite and the image of $F$ is closed, we deduce the same for $\mathfrak{F}$ 
using the orthogonality of the direct sum $U\oplus V$ and the fact that the 
image of $F^*$ is closed as well (see Theorem 3.2 in Lance \cite{Lance}).
Since $\mathfrak{F}$ satisfies the assumptions of Theorem 4, 
we have that $\mathfrak{F}$ is $A$-Fredholm. Consequently, there exists
a map $ \mathfrak{G}=\bigl(\begin{smallmatrix}
A & B \\ C& D
\end{smallmatrix} \bigr) \in \mbox{End}_{H_A^*}(U\oplus V) $ such that 
$\mathfrak{F}\mathfrak{G} = 1_{U\oplus V} + 
P_{U\oplus V}$ for an $A$-compact operator $P_{U\oplus V}$ in 
$U\oplus V.$ Expanding this equation, we get $FD = 1_V + 
\textrm{proj}_V {P_{U\oplus V}}_{|V}.$  
 It is immediate to realize 
that
$\textrm{proj}_V  {P_{U\oplus V}}_{|V}$ is an $A$-compact operator in $V.$ 
Thus,
$F$ is right invertible modulo an $A$-compact operator in $V.$ Similarly, one 
proceeds in the case of the left inverse. Summing-up, $F$ is an $A$-Fredholm 
morphism. 
\hfill $\Box$

\bigskip
In the next theorem, we study how certain  properties of the continuous 
extensions of 
pre-Hilbert module morphisms transfer to properties of the original map.
 
\bigskip
 
\noindent {\bf Theorem 6:} Let $A$ be a $C^*$-algebra of compact operators, $(V, 
(,)_{V})$ and 
$(W,(,)_{W})$  be Hilbert $A$-modules, and $(U,(,)_U)$ be a pre-Hilbert 
$A$-module 
which is a vector subspace of $V$ and $W$ such that the  norms 
 $|\,|_W$ and $|\,|_U$ coincide on $U$ and $|\,|_V$ restricted to $U$ dominates $|\,|_U.$
 Suppose 
that $D\in \textrm{End}_{PH_A^*}(U)$ is a self-adjoint 
morphism  
having a continuous adjointable extension $\widetilde{D}: V \to W$ such that 
\begin{itemize} 
\item[i)] $\widetilde{D}$ is $A$-Fredholm, 
\item[ii)] $\widetilde{D}^{-1}(U), \widetilde{D}^{*-1}(U) \subseteq U$ and
\item[iii)] $\textrm{Ker}\, \widetilde{D}$ and $\textrm{Ker}\, \widetilde{D}^*$ 
are 
subsets of $U.$
\end{itemize}
Then $D$ is a self-adjoint parametrix possessing operator in $U.$

{\it Proof.} 
We construct the parametrix and the projection.
\begin{itemize}

\item[1)] Using assumption (i), $\widetilde{D}$ has closed image by Corollary 
5. 
By 
Theorem 
3.2 
in Lance \cite{Lance}, 
the image of $\widetilde{D}^*:W \to V$ is closed as well, and  the following 
decompositions
$$V = \textrm{Ker} \, \widetilde{D}    \oplus  \textrm{Im} \, \widetilde{D}^*,$$
$$W = \textrm{Ker} \, \widetilde{D}^*   \oplus  \textrm{Im} \, \widetilde{D}  $$
hold.
Restricting $\widetilde{D}$ to the  Hilbert $A$-module $\textrm{Im}\, 
\widetilde{D}^*,$ we obtain a continuous bijective
Hilbert $A$-module homomorphism $\textrm{Im} \, \widetilde{D}^* \to 
\textrm{Im}\, 
\widetilde{D}.$

 Let us set
$$\widetilde{G}(x) = \left\{ \begin{array}{lll}  ({\widetilde{D}}_{| 
\textrm{Im} 
\, \widetilde{D}^*})^{-1}(x)&  \quad &x \in \textrm{Im}\ \widetilde{D}\\
                    0&   \quad &x \in \textrm{Ker}\, \widetilde{D}^*.
                   \end{array} \right. 
$$

The operator $\widetilde{G}:W \to V$ is continuous by  the  open map theorem. 
Due to its construction,
$\widetilde{G}$ is a morphism in the category $H_A.$  Because of the 
adjointability of
$\widetilde{D},$ and the definition of $\widetilde{G},$ $\widetilde{G}$ is 
adjointable as well. Summing-up, $\widetilde{G} \in \textrm{Hom}_{H_A^*}(W,V).$
Note that $\widetilde{G}: W \to \textrm{Im}\,\widetilde{D}^*.$

\item[2)] It is easy to see that the decomposition $V = \textrm{Ker} \, 
\widetilde{D} \oplus \textrm{Im} \, \widetilde{D}^*$
restricts to $U$ in the sense that $U = \textrm{Ker}\, D \oplus (\textrm{Im} \, 
\widetilde{D}^*\cap U).$ Indeed, let
$u\in U.$ Then $u \in V$ and thus $u = v_1 + v_2$ for $v_1 \in \textrm{Ker} \, 
\widetilde{D}$ and
$v_2 \in \textrm{Im} \, \widetilde{D}^*.$ Since $\textrm{Ker} \, \widetilde{D} 
\subseteq U$ (assumption (iii)) and $\textrm{Ker} \, D \subseteq \textrm{Ker}\, 
\widetilde{D},$  
we have $\textrm{Ker} \, \widetilde{D} = \textrm{Ker} \, D.$ 
Similarly, one proves that  $\textrm{Ker} \, \widetilde{D}^* = \textrm{Ker} \, 
D.$ In particular, $v_1 \in  \textrm{Ker}\, D.$ Since $U$ is a vector space,
$v_2 = u - v_1$ and $u, v_1 \in U,$  $v_2$ is an element of $U$ as well. Thus, 
$U \subseteq \textrm{Ker}\, D \oplus 
(\textrm{Im} \, \widetilde{D}^* \cap U).$ Since $\textrm{Ker}\, D, \textrm{Im} 
\, \widetilde{D}^* \cap U \subseteq U,$ the announced decomposition holds.

\item[3)]
Further, we have $\textrm{Im}\, \widetilde{D}^* \cap U = \textrm{Im}\, D.$ 
Indeed,
if $u\in U$ and $u = \widetilde{D}^*w $ for an element $w\in W$  then $w\in U$ 
due to item (ii) and consequently, $u=\widetilde{D}^*w = D^*w = Dw$ that 
implies $\textrm{Im}\, \widetilde{D}^* \cap U \subseteq \textrm{Im}\, D.$
The opposite inclusion is immediate.
(Similarly, one may prove that $\textrm{Im}\, \widetilde{D} \cap U = 
\textrm{Im}\, D.$)
Putting this result together with the conclusion of item 2 of this 
proof,  we obtain $U= \textrm{Ker} \, D \oplus \textrm{Im}\, D.$

\item[4)]
It is easy to realize that $\widetilde{G}_{|U}$ is into $U.$
Namely, if $v = \widetilde{G} u$ for an element $u \in U \subseteq V,$  we 
may write it as
$u = u_1 + u_2$ for $u_1 \in \textrm{Ker}\,  \widetilde{D}$ and $u_2 \in 
\textrm{Im}\, \widetilde{D}^*$ according to the decomposition of $V$ above.
Since $u_2 = u -u_1$ and $u_1 \in U$ (due to (iii)), we see that $u_2$ is an 
element of $U$ as well.
Consequently, $v = \widetilde{G}_{|U}u = 
\widetilde{G} u_1 + \widetilde{G} u_2 = \widetilde{D}^{-1}_{| \textrm{Im}\, 
{\widetilde{D}^*}}u_2.$ Since $\widetilde{D}^{-1}(U) \subseteq U$ (item (ii)), 
we 
obtain that $v\in U$ proving that $\widetilde{G}_{|U}$ is into $U.$ Let us set 
$G=\widetilde{G}_{|U}.$ Due to the assumptions on the norms and the continuity 
of $\widetilde{G}: (W,|\,|_W) \to (V,|\,|_V),$ it is easy to see that
$G:U\to U$ is continuous as well. 

\item[5)] Defining $P$ to be the projection of $U$ onto 
$\textrm{Ker}\, D$ along the $\textrm{Im}\, D,$   we get a self-adjoint  
projection on the pre-Hilbert module $U$ due to the decomposition $U= 
\textrm{Ker}\, D \oplus \textrm{Im}\, 
D$ derived in item 2 of this proof. The relations $DP=0$  and $1_U = GD + P = 
DG + P$ are then easily verified using the relation $\textrm{Ker} \, 
\widetilde{D}^* = \textrm{Ker} \, D.$

\end{itemize}
\hfill $\Box$

\bigskip

\noindent {\bf Remark 9:} In the preceding theorem,   
specific properties are generalized which are well  known to hold for self-adjoint elliptic operators 
acting on  smooth sections of vector bundles  over 
compact manifolds. For instance, assumption (ii) is the smooth regularity and (iii) expresses the fact that differential
 operators are of finite order. 
See, e.g., Palais \cite{Palais} or  Wells \cite{Wells}.

\section{Complexes of  pseudodifferential operators in $C^*$-Hilbert bundles}

For a definition of a $C^*$-Hilbert bundle, bundle atlae and  differential 
structures of bundles,
see Kr\'ysl \cite{KryslJGSP}, \cite{KryslAGAG2} or Mishchenko, Fomenko 
\cite{FM}.
For definitions of the other notions used in the next two paragraphs, we 
refer to  Solovyov, 
Troitsky \cite{ST}. 
Let us recall that for an $A$-pseudo\-diffe\-rential operator $D: 
\Gamma(M,\mathcal{E})\to \Gamma(M,\mathcal{F})$ acting between
 smooth  sections of 
$A$-Hilbert bundles $\mathcal{E}$ and $\mathcal{F}$ over a manifold $M,$ we 
have 
the order $\textrm{ord}(D) \in \mathbb{Z}$ of $D$
 and the symbol map $\sigma(D): \pi^*(\mathcal{E}) \to \pi^*(\mathcal{F})$ of 
$D$ 
at 
our disposal. Here, the map $\pi: T^*M \to M$ denotes the projection of the 
cotangent bundle. Moreover, if $M$ is compact, then for  $A$-Hilbert bundles 
$\mathcal{E} \to 
M$ and $\mathcal{F} \to M,$ an $A$-pseudodifferential operator 
$D:\Gamma(M,\mathcal{E}) 
\to \Gamma(M,\mathcal{F})$, 
and an integer $t\geq \textrm{ord}(D),$ we can form 
\begin{itemize}
\item[1)]the so called Sobolev type completions  $(W^t(M,\mathcal{E}),(,)_t)$   
of 
$(\Gamma(M,\mathcal{E}),(,))$  
\item[2)]the adjoint $D^*:\Gamma(M,\mathcal{F}) \to \Gamma(M,\mathcal{E})$ of 
$D$ 
and 
\item[3)] the continuous extensions
 $D_t:  W^t(M,\mathcal{E}) \to W^{t-\textrm{ord}(D)}(M,\mathcal{F})$ of $D.$
\end{itemize}

Smooth sections  $(\Gamma(M, \mathcal{G}),(,))$ of an 
$A$-Hilbert bundle $\mathcal{G} \to M$ form a pre-Hilbert $A$-module and
  spaces $(W^t(\mathcal{G}),(,)_t)$ are Hilbert $A$-modules. See Example 2 
item 4 for a definition of the $A$-product $(,)$ on the space of smooth 
sections.
The adjoint $D^*$ of an $A$-pseudodifferential operator $D$ 
is considered with respect to the $A$-products $(,)$ on the pre-Hilbert 
$A$-modules of smooth sections of the appropriate bundles. 
  Operators  $D$ and $D^*$ are pre-Hilbert $A$-module  morphisms,  
extensions $D_t$ are Hilbert $A$-module  morphisms, 
and the symbol map $\sigma(D)$ is a morphism of $A$-Hilbert bundles.

\bigskip

The $A$-ellip\-ticity is defined similarly as the ellip\-ticity of differential 
operators in  bundles 
with 
finite dimensional fibers  over $\mathbb{R}$ or $\mathbb{C}.$
We use the following definition, the first part of which is contained  in 
Solovyov, Troitsky \cite{ST}. 

\bigskip

\noindent {\bf Definition 7:} Let $D:\Gamma(M,\mathcal{E}) \to 
\Gamma(M,\mathcal{F})$ be 
an 
$A$-pseudodifferential
operator. We say that $D$  is $A$-{\it elliptic} if $\sigma(D)(\xi, 
-):\mathcal{E} \to \mathcal{F}$ is an isomorphism of $A$-Hilbert bundles for 
any 
non-zero
$\xi \in T^*M.$ Let $(p_i:\mathcal{E}^i \to M)_{i\in \mathbb{Z}}$ be a 
sequence of $A$-Hilbert bundles and
$(\Gamma(M,\mathcal{E}^i), d_i: \Gamma(M,\mathcal{E}^i) \to 
\Gamma(M,\mathcal{E}^{i+1}))_{i\in \mathbb{Z}}$ be a complex of 
$A$-pseudodifferential operators.
We say that  $d^{\bullet}$ is {\it $A$-elliptic}  if and only if the complex of 
symbol maps $(\mathcal{E}^i,\sigma(d_i)(\xi,-))_{i\in \mathbb{Z}}$
is exact for each non-zero $\xi \in T^*M.$

\bigskip

\noindent {\bf Remark 8:} One can show that  the Laplace operators 
$\Delta_i = d_{i-1}d_{i-1}^* + d_i^*d_i,$ $i \in \mathbb{Z},$ of an 
$A$-elliptic complex are $A$-elliptic
operators in the sense of Definition 7. For a proof in the $C^*$-case, see 
Lemma 9 in Kr\'ysl \cite{KryslAGAG1}. Let us notice that the assumption on 
unitality of $A$ is inessential in the  proof of the Lemma 9 in 
\cite{KryslAGAG1}.

\bigskip

Recall that an $A$-Hilbert bundle $\mathcal{G} \to M$ is called {\it 
finitely generated 
projective} if its fibers are finitely generated and projective
Hilbert $A$-modules. See  Manuilov, Troitsky \cite{MT}.
Let us recall a theorem of Fomenko and Mishchenko on a relation of 
the $A$-ellip\-ticity and the $A$-Fredholm property.

\bigskip

\noindent {\bf Theorem 7:} Let $A$ be a $C^*$-algebra, $M$ a compact manifold, 
$\mathcal{E} \to M$ a finitely generated
projective $A$-Hilbert bundle over $M,$ and $D: \Gamma(M,\mathcal{E}) \to 
\Gamma(M,\mathcal{E})$ an $A$-elliptic
operator. Then the continuous extension $$D_t:W^t(M,\mathcal{E}) \to W^{t - 
\textrm{ord}(D)}(M,\mathcal{E})$$ is an $A$-Fredholm morphism for any 
$t \geq \textrm{ord}(D).$ 

{\it Proof.} See Fomenko, Mishchenko \cite{FM} and Remark 6.
\hfill $\Box$

\bigskip

\noindent {\bf Corollary 8:} Under the assumptions of Theorem 7, 
$\textrm{Ker}\, D_t = 
\textrm{Ker}\, D$ for any $t\geq \textrm{ord}(D).$
If moreover $D$ is self-adjoint, then also $\textrm{Ker} \, {D_t}^* = 
 \textrm{Ker}\, D$ for any $t \geq \textrm{ord}(D).$

{\it Proof.} See Theorem 7 in Kr\'ysl \cite{KryslAGAG1} for the first claim, 
and the formula (5)  in \cite{KryslAGAG1} for the second one. 
\hfill $\Box$

\bigskip

Let us notice that the first assertion in Corollary 8 appears as Theorem 
3.1.145 on  pp. 101 in Solovyov,  Troitsky \cite{ST}.
Now, we state the theorem saying that the Hodge theory holds for
$A$-elliptic complexes of operators acting on sections of finitely generated 
projective $C^*$-Hilbert bundles over compact manifolds if $A$ is a 
$C^*$-algebra of compact operators.

\bigskip

\noindent {\bf Theorem 9:} Let $A$ be a $C^*$-algebra of compact operators,
$M$ be a compact manifold, $(p_i:\mathcal{E}^i \to M)_{i\in \mathbb{Z}}$ be a 
sequence of finitely generated projective $A$-Hilbert bundles over $M$ and 
$d^{\bullet}= (\Gamma(M,\mathcal{E}^i), d_i: \Gamma(M,\mathcal{E}^i) \to 
\Gamma(M,\mathcal{E}^{i+1}))_{i\in \mathbb{Z}}$ be a complex of 
$A$-pseudodifferential operators.
If  $d^{\bullet}$ is $A$-elliptic,
then  for each $i\in \mathbb{Z}$
\begin{itemize}
\item[1)] $d^{\bullet}$ is of Hodge type, i.e., $\Gamma(M,\mathcal{E}^i) = 
\textrm{Ker} \, \Delta_i \oplus 
\textrm{Im}\, d_i^* \oplus \textrm{Im} \, d_{i-1}$
\item[2)] $\textrm{Ker}\, d_{i} = \textrm{Ker} \, \Delta_i \oplus 
\textrm{Im} 
\, d_{i-1}$
\item[3)] $\textrm{Ker}\, d_{i}^* = \textrm{Ker} \, \Delta_{i+1} \oplus 
\textrm{Im} \, d_{i+1}^*$
\item[4)] $\textrm{Im}\, \Delta_i  = \textrm{Im} \, d_{i-1} \oplus 
\textrm{Im} \,  d_i^*$
\item[5)] The cohomology group $H^i(d^{\bullet})$ is a finitely generated 
projective 
$A$-Hilbert module isomorphic to 
the $A$-Hilbert module $\textrm{Ker} \, \Delta_i.$
\end{itemize}
 
{\it Proof.} Since $d^{\bullet}$ is an $A$-elliptic complex, the associated  
Laplace operators are $A$-elliptic operators (Remark 8). 
The Laplace operators are self-adjoint according to their definition. 
According to Theorem 7, the extensions $(\Delta_i)_t$ are 
$A$-Fredholm for any $t\geq \textrm{ord}(\Delta_i).$

   Let us set $D=\Delta_i,$ $U=\Gamma(M,\mathcal{E}^i),$ $V = 
W^{\textrm{ord}(\Delta_i)}(M, \mathcal{E}^i)$ and 
$W=W^0(M, \mathcal{E}^i)$  considered with the appropriate $A$-products. 
Then   $U$ is a vector subspace of $V \cap W,$ and the 
restriction of $(,)_W$ to $U\times U$ coincides with
$(,)_U.$
Since $\Delta_i$ is $A$-elliptic, $\textrm{Ker} \, \Delta_i = \textrm{Ker}\, 
(\Delta_i)_t = \textrm{Ker}\, 
(\Delta_i)_t^*$ due to Corollary 8.
Because the operator  $D$ is of finite order, 
$\widetilde{D}^{-1}(\Gamma(M,\mathcal{E}^i)),$ 
$\widetilde{D}^{*-1}(\Gamma(M,\mathcal{E}^i))  \subseteq  
\Gamma(M,\mathcal{E}^i).$
   The norm on $U=\Gamma(M,\mathcal{E}^i)$ coincides with the norm on 
$W=W^0(M,\mathcal{E}^i)$ restricted to $U$ and the norm $|\,|_U$ on $U$ is 
dominated by the 
norm $|\,|_V$ on $V=W^{\textrm{ord}(D)}(M,\mathcal{E}^i)$ restricted to $U.$  
Thus, the assumptions on the norms in Theorem 6 are satisfied
and we may conclude, that 
$\Delta_i$ is a self-adjoint parametrix possessing morphism, and thus,
$d^{\bullet}$ is self-adjoint parametrix possessing as well.

   The assertions in items 1--4  follow from the corresponding assertions of  
Theorem 2. Using Theorem 2 item 5, $H^{i}(d^{\bullet}) \simeq \textrm{Ker} \, 
\Delta_i.$ As already mentioned, $\textrm{Ker}\, \Delta_i \simeq \textrm{Ker}\, 
(\Delta_i)_t.$ Since $(\Delta_i)_t$ is $A$-Fredholm and $A$ is a $C^*$-algebra 
of compact operators, $\textrm{dim}_A \textrm{Ker} (\Delta_i)_t$ is finite due 
to Corollary 5. It follows that the kernel of $\Delta_i$ is finitely generated.

   Since the image of $(\Delta_i)_t$ is closed (Corollary 5), we have 
$W^t(M,\mathcal{E}^i) 
= \textrm{Ker}\, (\Delta_i)_t \oplus \textrm{Im} \, 
(\Delta_i)_t^*$ due to the Mishchenko theorem (Theorem 3.2 in Lance 
\cite{Lance}).  Consequently,  $\textrm{Ker}\, \Delta_i = 
\textrm{Ker}\, (\Delta_i)_t$ is a  projective $A$-Hilbert module by Theorem 
1.3 in Fomenko, 
Mishchenko \cite{FM}.
\hfill$\Box$

\bigskip

\noindent {\bf Remark 10:} Let us notice that if the assumptions  of Theorem 9 
are satisfied, the cohomology groups  share the properties of their fibers 
in the sense that they are   finitely generated projective 
$A$-Hilbert modules.

In the proof of Theorem 9, we could have shown  the cohomology groups to be
finitely generated and projective in a shorter way using  Theorem 7 with
an alternative definition of an $A$-Fredholm operator given in 
\cite{FM} since also Theorem 7 as appears in \cite{FM} uses the 
alternative notion of an $A$-Fredholm  map.

\bibliographystyle{amsplain}

\end{document}